\documentclass[]{article}

\usepackage{amsfonts,amsmath,amssymb}
\usepackage{amsthm}
\usepackage[cp866]{inputenc}
\usepackage[english]{babel}

\usepackage{xcolor}

\allowdisplaybreaks \theoremstyle{definition}

\newtheorem{Remark}{Remark}

\theoremstyle{plain}
\newtheorem{Theorem}{Theorem}

\newtheorem{Proposition}{Proposition}

\newtheorem{Fact}{Fact}

\newcommand \rat {{\mathbb Q}}

\newcommand \ls {\leqslant}

\newcommand \bx {\bar {x}}

\newcommand \ca {{\mathbb C}_{\text {\rm alg}}}

\newcommand \caa {C_{\text {\rm alg}}}

\title{Bit Complexity  of Computing  Solutions \\ for Symmetric Hyperbolic Systems of PDEs \\ with Guaranteed Precision}
\author{Svetlana Selivanova
\thanks{Partially supported  by RFBR-JSPS Grant 20-51-50001, by the National Research Foundation of
Korea (grant 2017R1E1A1A03071032), by the NRF Brain Pool program (grant 2019H1D3A2A02102240), by the  International Research \& Development Program of the Korean Ministry of Science and ICT (grant 2016K1A3A7A03950702), and by RFBR grant 17-01-00801.}
 \\ Korean Advanced Institute of Science and Technology, Republic of Korea\\
   {\tt sweseliv@gmail.com}
\\ and
\\Victor  Selivanov
\thanks{Partially supported by RFBR-JSPS Grant 20-51-50001 and by the Regional Mathematical Center of Kazan Federal University.}
 \\A.P. Ershov Institute of
Informatics Systems SB RAS, Russia
\\{\tt vseliv@iis.nsk.su}
}

\begin{document}
\large
\date{}
 \maketitle

%------------------

\begin{abstract}
We establish upper bounds on bit complexity  of computing solution operators for symmetric hyperbolic systems of PDEs, combining symbolic and approximate algorithms to obtain the solutions with guaranteed prescribed precision. Restricting to algebraic real inputs allows us to use the classical (``discrete'') bit complexity concept.


%{\bf MSC}: Primary   03D15, 03D78, 58J45; Secondary 65M06, 65M25.
\end{abstract}

\section{Introduction}
In this paper, which extends the conference paper \cite{ss18}, we establish upper bounds on (classical ``discrete'') bit complexity  of computing solution operators for symmetric hyperbolic systems of PDEs (restricting the consideration to algebraic real inputs). 

Algorithms used in mathematics-oriented software can be divided into two big classes:  symbolic algorithms which aim to find  precise solutions, and  approximate algorithms which aim to find  ``good enough'' approximations to precise solutions. The symbolic algorithms are implemented e.g. in computer algebra systems while the approximate algorithms are included into numerical mathematics packages. Both classes of algorithms are widely used in applications and in mathematical research. 

Symbolic algorithms correspond well to computations on discrete structures (with mathematical foundations in the classical computability and complexity theory) while the approximate algorithms help to carry out computations on continuous structures (with mathematical foundations in the field of computability and complexity in analysis which however is not yet completely bridged with numerical methods used in practical applications).
An important idea relating  both classes of algorithms is to look for  approximations to the precise solutions  with ``guaranteed precision''. It is one of the approaches in the fast developing branch of reliable computations which are crucially important for safety-critical applications. 

The statement of a guaranteed-precision version of some problem on a continuous structure may  reduce it to a problem on a discrete structure which enables to apply the classical computability and bit complexity theory. 
This is the approach 
we adopt here to measure the complexity of algorithms for differential equation problems: searching for an approximate solution with guaranteed prescribed precision and using exact symbolic calculations on each step of the algorithm. Classical bit complexity is fundamental because it estimates the amount of computational resources needed to solve a problem on a computing device. We briefly recall the definitions in Section \ref{coding} (also there are comments on other complexity approaches) and explain how to encode algebraic reals: using encodings of the corresponding minimal (rational) polynomial and the number of the root to which the given algebraic real corresponds, see Section \ref{algprel}.

In this paper,  we 
consider initial-value (IVP), or Cauchy, and well posed boundary-value (BVP) problems for symmetric hyperbolic systems of PDEs $A\frac{\partial
{\bf u}}{\partial t}+\sum\limits_{i=1}^mB_i\frac{\partial {\bf u}}
{\partial x_i}=f(t,x_1,\ldots,x_m)$ where $A=A^\ast>0$ and $B_i=B_i^\ast$ are symmetric
 $n\times n$-matrices, $t\geq0$, $x=(x_1,\ldots,x_m)\in
Q=[0,1]^m$,  $f:[0,+\infty)\times Q\rightharpoonup{\mathbb R}^n$. Such systems   can be used to describe a wide variety of physical
processes like those considered in the theories of elasticity,
acoustics, electromagnetism etc., see e.g.
\cite{fri,go71,go76,kps}. Using a well known grid-based method (difference scheme), for algebraic real coefficients and rational polynomials as initial data, we establish an EXPTIME complexity upper bound, and indicate restrictions when it boils down to PTIME:  the precision needs to be fixed and certain upper bounds on inputs need to be imposed.

We consider both the IVP and BVP
cases as very general and illustrative: they are different in terms of domains of existence and uniqueness as well as properties of the solutions, but 
 almost the same w.r.t. the difference approach. The methods of the paper are applicable to IVP and BVP  not only for symmetric hyperbolic systems, but for broader classes of PDE systems. 
 
 We stick here to symmetric hyperbolic systems as a quite broad and practically important class capturing one of the following common problems with (explicit) difference scheme for evolutionary systems of PDEs. 
 To construct a stable difference scheme it is necessary to know the number of positive/ negative/ zero eigenvalues and corresponding eigenvectors of matrix pencils related to $A,B_i$, which are known to be non-computable unless e.g. the number of different eigenvalues is given as input \cite{zb04}. The computational instabilities occuring because of that are well known in practice; here we are carefully avoiding them by restricting to algebraic real inputs, while in \cite{ss13,ss17} we add extra restrictions for rigorously proving computability of solution operators in the real setting, see Subsection \ref{cbv}.

Moreover, for the considered case of symmetric hyperbolic systems we provide an explicit PTIME algorithm of computing, from the given precision and input data, the (space and time) grid steps, using the difference scheme with which (or any smaller) provides the solution with this given precision, see Proposition \ref{steps_prop} in Subsection \ref{steps}. A similar approach can be used for broader classes of systems of PDEs, as well as for the real setting \cite{ss13}, which helps to adapt the difference scheme approach to exact real computation (ERC) packages, see \cite{sz}. In fact, while the present paper was being reviewed, the first author has published, with two other coauthors, the paper 
\cite{ksz}
with real complexity estimates for linear evolutionary systems of PDEs (including symmetric hyperbolic systems as a particular case) in the exact real computation approach. The therein obtained estimates help to improve the ones in the present paper: PSPACE instead of EXPTIME, due to replacing the step-by-step difference scheme iterations by a more efficient procedure (recursive matrix powering). However, in \cite{ksz} we consider fixed real-valued initial functions and matrix coefficients while in the present paper we treat them as inputs, so these papers complement each other; in the difference scheme part \cite{ksz} partially uses the methods of the present paper.

In the next section we formally state the problems to be investigated, and recall some notions and background to be used. In Section \ref{aux} we establish  upper complexity bounds for several auxiliary linear algebra problems closely related to finding the solutions of symmetric hyperbolic systems of PDEs. In Section \ref{mainres} we prove the bit cost estimates for the solutions of PDEs, and we conclude in Section \ref{con} with a short discussion on possible applications of this work.

\section{Preliminaries and formulations of main results}\label{prel}

\subsection{Cauchy and boundary-value problems}\label{cbv}

The Cauchy problem for a symmetric hyperbolic system is stated as follows:
\begin{equation} \label{sist_1}
\begin{cases} A\frac{\partial
{\bf u}}{\partial t}+\sum\limits_{i=1}^mB_i\frac{\partial {\bf u}}
{\partial x_i}=f(t,x_1,\ldots,x_m),\ t\geq 0,\\
{\bf u}|_{t=0}=\varphi(x_1,\ldots,x_m),
\end{cases}
\end{equation}
where $A=A^\ast>0$ and $B_i=B_i^\ast$ are  symmetric
 $n\times n$-matrices, $t\geq0$, $x=(x_1,\ldots,x_m)\in
Q=[0,1]^m$, $\varphi:Q\rightarrow{\mathbb R}^n$, $f:[0,+\infty)\times Q\rightharpoonup{\mathbb R}^n$ and ${\bf
u}:[0,+\infty)\times Q\rightharpoonup{\mathbb R}^n$ is a partial
function acting on the domain $H$ of existence and uniqueness of
the Cauchy problem \eqref{sist_1}.  The set $H$ is known to be (see e.g.
\cite{go71}) the intersection of  semi-spaces
 $$
t\geq0,\;x_i-\mu^{(i)}_{\rm
max}t\geq0,\;x_i-1-\mu^{(i)}_{\rm min}t\leq0\;(i=1,\ldots,m)
 $$
of ${\mathbb R}^{m+1}$ where $\mu^{(i)}_{\rm
min},\mu^{(i)}_{\rm
max}$ are respectively the minimum and maximum of the eigenvalues of the matrix  $A^{-1}B_i$.

The boundary-value problem is stated as follows:
\begin{equation}\label{sist_2}
\begin{cases} A\frac{\partial
{\bf u}}{\partial t}+\sum\limits_{i=1}^mB_i\frac{\partial {\bf u}}
{\partial x_i}=f(t,x_1,\ldots,x_m), \\
{\bf u}|_{t=0}=\varphi(x_1,\ldots,x_m), \\
\Phi_i^{(1)}{\bf u}(t,x_1,\ldots,x_{i-1},0,x_{i+1},\ldots,x_m)=0,\\
\Phi_i^{(2)}{\bf u}(t,x_1,\ldots,x_{i-1},1,x_{i+1},\ldots,x_m)=0,\\
i=1,2,\ldots,m,
\end{cases}
\end{equation}
where $A,B_i,Q,\varphi,f,\mathbf{u}$ are as above, and $\Phi_i^{(1)}$, $\Phi_i^{(2)}$ are rectangular matrices (the boundary coefficients) that  meet the
following conditions:

1) The number of rows of $\Phi_i^{(1)}$ (respectively,
$\Phi_i^{(2)}$) is  equal to the number of positive (respectively,
negative) eigenvalues of the matrices $A^{-1}B_i$; coincidence constraints of the initial and boundary conditions hold (such constraints depend on the particular problem and on the smoothness which we want to obtain).

2) The boundary conditions are assumed to be dissipative which
means that
\begin{equation}\label{gran_dissip}
(B_i{\bf u},{\bf u})\leq 0\text{ for }x_i=0,\quad (B_i{\bf u},{\bf
u}) \geq 0\text{ for }x_i=1,\quad i=1,2,\ldots,m.
 \end{equation}
 
Condition 1) guarantees the existence of solution of the boundary-value problem \eqref{sist_2} in the cylinder $[0,\infty)\times Q$, while 2) implies its uniqueness. 

For both problems \eqref{sist_1} and \eqref{sist_2}, theorems on continuous dependence of  solutions on the input data hold, i.e. the problems are correctly posed. Both problems are practically important, since many physical processes (including linear elasticity, acoustic, Maxwell equations) are described by such kind of systems.

\medskip
For the considered problems, there exist different numerical methods, from which we use those developed in \cite{go76,go71}, see also \cite{kps}. Their convergence relies on the following well-known theorem (see e.g. \cite{gr,kps,strik}):  if a difference scheme approximates the given differential problem and is stable (which is an intrinsic property of the difference scheme), then the (discrete) solution of the corresponding difference equations converges to the exact solution of the differential problem in an appropriate grid norm; the speed of convergence corresponds to the order of approximation. 

The Godunov scheme, which we use, is of first order of approximation and is stable with a Courant number (relating the time and space steps of the grid) depending on spectral characteristics of the matrix coefficients $A$, $B_i$ (see Section \ref{steps} below). As it is known, in order to construct a stable difference scheme for a symmetric hyperbolic system, one needs to compute eigenvectors of symmetric matrices, which is actually a discontinuous operation \cite{rel}, hence not computable. However, due to \cite{zb04,zie}, eigenvectors are computable, provided that the spectrum cardinality  (i.e., the number of distinct eigenvalues of the given matrix) is given as an input.

In \cite{ss09, ss13, ss17} we developed an approach to study 
computability properties of PDEs based on the Godunov difference scheme and established the
computability, in the sense of the TTE approach \cite{bhw,wei}, of the solution
operator  of the Cauchy and boundary-value problems. In this section we briefly recall main properties of such systems and the related computability results.

For computability of the solution of \eqref{sist_1} from the initial data and the matrix coefficients, the following result has been established.

\begin{Theorem}\cite{ss13}\label{main2}
Let $M_{\varphi}>0, M_A>0, p\geq2$ be integers,   let
$i=1,\ldots,m$, and let $n_A,n_1,\ldots,n_m$ be cardinalities of spectra of $A$ and of the matrix pencils $\lambda A- B_1,\ldots,\lambda A-B_m$, respectively (i.e., $n_i$ is the number of distinct roots of the characteristic polynomial
${\rm det}(\lambda A- B_i)$). Then the operator
 $$R: (A,B_1,\ldots,B_m,n_A,n_1,\ldots,n_m,\varphi)\mapsto{\mathbf u}$$
  sending any  sequence
$A,B_1,\ldots,B_m$ of symmetric real matrices with $A> 0$ such that the matrix pencils $\lambda A- B_i$ have no zero eigenvalues,
 \begin{equation}\label{M_A}
 ||A||_2,||A^{-1}||_2,||B_i||_2\leq M_A,\quad
\lambda^{(i)}_{\min}<0<\lambda^{(i)}_{\max},\quad i=1,2,\ldots,m,
\end{equation}
 the sequence $n_A,n_1,\ldots,n_m$ of the corresponding cardinalities, and any function $\varphi\in C^{p+1}(Q,{\mathbb R}^n)$ satisfying
the conditions 
\begin{equation}\label{bound}
||\frac{\partial
\varphi}{\partial x_i}||_{s}\leq M_{\varphi},\ ||\frac{\partial^2
\varphi}{\partial x_i\partial x_j}||_{s}\leq M_{\varphi},\
 i,j=1,2,\ldots,m,
 \end{equation} to the unique solution ${\mathbf
u}\in C^p(H,{\mathbb R}^n)$ of \eqref{sist_1} (with $f=0$) is a computable
partial function from the space $S_{+}\times S^{m}\times \mathbb{N}^{m+1}\times
C_s^{p+1}(Q,{\mathbb R}^n)$ to $C_{sL_2}^p(H,{\mathbb R}^n)$.
\end{Theorem}

Here  $C_s^{p+1}(Q,{\mathbb R}^n)$, $C^p_{sL_2}(H,{\mathbb R}^n)$ are the spaces of continuously differentiable functions, being in $C^{p+1}$ or $C^p$ on the corresponding sets, such that all of their first and second derivatives are uniformly bounded by $M_{\varphi}$, with  the sup-norm  $$||\varphi||_s=\sup_{x\in
Q}|\varphi(x)|$$ on
$C^{p+1}(Q,\mathbb R^n)$ and 
$sL_2$-norm
 $$||{\bf u}||_{sL_2}=\sup_{\{t:(t,x)\in H\}}\sqrt{\int_Q|{\bf u}(x,t)|^2dx}$$
on $C^p(H,\mathbb R^n)$, respectively;
$|{\bf u}(x,t)|^2=\langle{\bf u},{\bf u}\rangle$ is the standard scalar product.

 By $S$ and $ S^{+}$ we denote respectively the spaces of all symmetric and symmetric positively definite matrices with euclidean norms uniformly bounded by a constant $M_A$.

In \cite{ss13}, an analogue of Theorem \ref{main2}  is established for the boundary-value problem \eqref{sist_2} with  fixed
computable real matrices $\Phi_i^{(1)}$, $\Phi_i^{(2)}$ $(i=1,2,\ldots,m)$,  provided that the strong dissipativity property (with strict inequalities in \eqref{gran_dissip}) holds. The solution for the boundary-value problem is defined on a cylinder $[0,T]\times Q$, for any computable $T>0$, i.e.  ${\bf u}\in C^p_{sL_2}([0,T]\times Q,\mathbb
R^n)$.

\begin{Remark}
 The computability of the operator $R: (A,B_1,\ldots,B_m,\varphi)\mapsto{\mathbf u}$ without the additional information on spectra cardinalities is by now an open question in general, though in particular cases like Cauchy or periodic boundary-value problems computability can be established by different methods like e.g. Fourier transformation. In case of non-computability, it is instructive to study the topological complexity of computing the operator.
For the case when the number of different eigenvalues of a matrix is given as input, a PTIME algorithm of computing eigenvectors has been recently devised by Sewon Park (the paper soon to appear); also eigenvalues are known to be PTIME computable  \cite{yap}, see also \cite{panchen}. Combined with proofs of our Theorems \ref{main2}, \ref{cpq} these facts give an analog of Theorem \ref{cpq} in the real complexity approach; those estimates can then be further improved as in \cite{ksz}. For different approaches to define complexity  see Remark \ref{complex_remark}.
\end{Remark}

Interestingly, if we work in decidable (in the Russian terminology, strongly constructivizable) fields (for example, the field of algebraic reals $\mathbb A$ used in the present paper), we do not need to add cardinalities as inputs.

\begin{Theorem}\cite{ss13}\label{commain}
Let $(\mathbb{B},\beta)$ be a decidable (i.e., strongly  constructive) real
closed ordered subfield  of $\mathbb{R}$. Then the solution operator $R: (A,B_1,\ldots,B_m,\varphi)\mapsto${\bf u} is uniformly computable (w.r.t. the numbering $\beta$) on  matrices $A,
B_1,\ldots,B_m$ with coefficients in $\mathbb{B}$. 
 \end{Theorem}

For the boundary-value problem  \eqref{sist_2} an analogue of this result is proved, with uniformity on the matrices $\Phi_i^{(1)}$, $\Phi_i^{(2)}$ having coefficients in $\mathbb{B}$. 

In \cite{ss17}, dependence on the right-hand part $f$ in \eqref{sist_1} is added, and also the time $T>0$ is included into the arguments of the solution operator $R$, i. e. the solution belongs to  
$C^p_{sL_2}(H\cap [0,T]\times Q,{\mathbb
R}^n)$. In \cite{ss17} the technique of using decidable (i.e., strongly constructive) fields for the computability results is explained in detail.

Also, Theorem \ref{commain} implies computability for fixed computable real matrices (which we earlier established using the result of \cite{zb04}).

\begin{Theorem}\cite{ss09,ss13}\label{main}
Let $Q=[0,1]^m$; $T>0$ be a computable real and $M_{\varphi}>0$, $p\geq 2$ be integers.
Let  $A, B_1,\ldots,B_m$ be fixed computable symmetric
matrices, such that $A=A^\ast> 0$, $B_i=B_i^\ast$.
 
 If $\varphi\in C^{p+1}(Q)$ satisfies
 the conditions \eqref{bound},
  then the operator
$R:\varphi \mapsto {\bf u}$ mapping the initial function $\varphi$ to  the unique solution ${\mathbf
u}\in C^p(H,{\mathbb R}^n)$ of the
Cauchy problem \eqref{sist_1} is a computable partial
function from $C^{p+1}_s(Q,\mathbb R^n)$ to $C^p_{sL_2}(H,\mathbb
R^n)$.
\end{Theorem}

An analogue of Theorem \ref{main}  is established in \cite{ss13} for the boundary-value problem \eqref{sist_2} with  fixed
computable real matrices $\Phi_i^{(1)}$, $\Phi_i^{(2)}$ $(i=1,2,\ldots,m)$.

The next natural step is to study computational complexity of  the mentioned problems, which we start to do in the present paper.

\subsection{Discretization of the problems}\label{discr}

To investigate the complexity of computing solutions of the systems \eqref{sist_1} and \eqref{sist_2} (and even to formulate the results), we  need discrete approximations of the given and unknown functions, as well as their interpolations. Therefore we
 first recall some discretization details. 
 
Consider, for any positive integer $N$, the uniform rectangular  grid $G_N$ on $Q=[0,1]^m$ defined
by the  points $$\left(\frac{i_1-\frac{1}{2}}{2^N},\frac{i_2-\frac{1}{2}}{2^N},\ldots,\frac{i_m-\frac{1}{2}}{2^N}\right)$$ where $1\leq
i_1,i_2, \ldots, i_m\leq2^N$. Let
$h=1/2^N$ be the corresponding spatial grid step and  $\tau$ be a time step. Denote  $G_N^{\tau}=G_N\times\{l\tau\}_{l=1}^L$, where $L$ is the number of the time steps.  
The choice of steps $h$ and $\tau$, which guarantee good properties of the difference scheme, will be specified below in Section \ref{steps}.

Note that the number of points in the grid $G_N$
is $2^{Nm}$, so the set $\mathbb{Q}^{G_N}$ (or $\mathbb{A}^{G_N}$) of grid functions
$g^{(h)}:G_N\to\mathbb{Q}^n$  (resp. $g^{(h)}:G_N\to\mathbb{A}^n$) may be identified with $\mathbb{Q}^{n\cdot
2^{Nm}}$ (resp. $\mathbb{A}^{n\cdot 2^{Nm}}$). We will consider  the following grid norms
 $$||g^{(h)}||_s=\operatorname{max}_{x\in G_N}|g^{(h)}(x)|, \;
 ||g^{(h)}||^2_{L_2}=h^m\sum_{x\in G_N}\langle g^{(h)}(x),g^{(h)}(x)\rangle.$$

We will consider the $sL_2$-norm on the vector spaces  $\mathbb{Q}^{G^\tau_N}$ (or $\mathbb{A}^{G^\tau_N}$ ) of grid
functions $v^{(h)}(t,x)$ on such grids: 
 $$||v^{(h)}||_{sL_2}=\operatorname{max}_{t\in\{l\tau\}_{l=1}^M}h^m\sum_{x\in G_N}\langle v^{(h)}(t,x),v^{(h)}(t,x)\rangle.$$

Recall that multilinear interpolations $\tilde{{\bf u}}$ (linear on each coordinate and coinciding with ${\bf u}$ at the grid points) have the following properties: ${\bf u}\mapsto\tilde{{\bf u}}$ and
${\bf u}^{(h)}\mapsto\widetilde{{\bf u}^{(h)}}$ are linear, and the following estimate holds \cite{za,sz59}: 
 \begin{equation}\label{int-est}||{\bf u}-\tilde{{\bf u}}||_s\leq \max_{i,j}\left\{||\frac{\partial^2{\bf u}}{\partial x_i\partial x_j}||_{sL_2}, ||\frac{\partial^2{\bf u}}{\partial x_i\partial t}||_{sL_2}\right\}\cdot h^2.\end{equation}

Further in Section \ref{mainres} we will construct, by means of a stable difference scheme approximating the differential system \eqref{sist_1},   a grid function $\upsilon$ on $G^{\tau}_N$ such that 
\begin{equation}\label{est}||{\bf u}-\widetilde{\upsilon\mid_H}||_{sL_2}<\frac{1}{a},\end{equation} where $a>1$ is a given integer (which determines the precision $\frac{1}{a}$  of computation) and ${\bf u}$ is the solution of \eqref{sist_1}.

Recall that, for an abstract boundary-value problem \begin{equation} \label{sist}
\begin{cases} L
{\bf u}(y)=f(y)\in C^p(\Omega,\mathbb R^n),\quad y\in\Omega\subset\mathbb R^{k}\\
{\mathcal L}{\bf u}(y)|_{\Gamma}=\varphi(y\mid_{\Gamma})\in
C^q(\Gamma,\mathbb R^n),\ \Gamma\subseteq \partial\Omega.
\end{cases}
\end{equation}
(where $L$ and ${\mathcal L}$ are differential operators with the
differential order of ${\mathcal L}$  less than that of $L$, $\Gamma$ is a part of the boundary $
\partial\Omega$ of some area $\Omega$), a {\it difference scheme} is a system of  algebraic equations 
\begin{equation}\label{sist_dif}L_h{\bf u}^{(h)}={\bf f}^{(h)},\;{\mathcal L}_h{\bf
u}^{(h)}=\varphi^{(h)}.
\end{equation}

Here $L_h,{\mathcal L}_h$ are difference operators (which are in
our case linear); $u^{(h)}$ and $\varphi^{(h)}$ are grid functions.

The scheme \eqref{sist_dif} {\it approximates} the differential equations \eqref{sist}
with order of accuracy $l$ (where $l$ is a positive integer) on a
solution ${\bf u}(t,x)$ of (\ref{sist}) if
\begin{eqnarray*}
 ||(L{\bf
u})|_{G_k}-L_h{\bf u}^{(h)}||_{F_h}\leq M_1h^l,\;
||f|_{G_k}-f^{(h)}||_{F_h}\leq M_2h^l,\\ ||({\mathcal L}{\bf
u})|_{G_k}-{\mathcal L}_h{\bf u}^{(h)}||_{\Phi_h}\leq M_3h^l \text{
and
 }||\varphi|_{G_k}-\varphi^{(h)}||_{\Phi_h}\leq M_4h^l
\end{eqnarray*}
 for some
constants $M_1,M_2,M_3$ and $M_4$ not depending on $h$ and
$\tau$.

The difference scheme \eqref{sist_dif} is called {\em stable} if its
solution ${\bf u}^{(h)}$ satisfies
$$||{\bf u}^{(h)}||_{U_h}\leq
N_1||f^{(h)}||_{F_h}+N_2||\varphi^{(h)}||_{\Phi_h}$$
for some
constants $N_1$ and $N_2$ not depending on $h$, $\tau$, $f^{(h)}$
and $\varphi^{(h)}$.

\begin{Fact}\cite{gr}\label{scheme} Let the difference scheme  be stable and
approximate (1) on the solution ${\bf u}$ with order $l$.
Then the solution $u^{(h)}$ uniformly converges to the
solution ${\bf u}$  in the sense that $||{\bf u}|_{G^\tau_k}-{\bf
u}^{(h)}||_{U_h}\leq Nh^l$ for some constant $N$ not depending on
$h$ and $\tau$.
\end{Fact}

The difference scheme which we use for approximating \eqref{sist_1} will be described in detail in Section \ref{mainres}.

\subsection{Algebraic preliminaries}\label{algprel}

In the study of computability of solution operators we considered rather general classes of initial data $A,B_i,\varphi_i,\cdots$ (matrices with real coefficients, broad classes of smooth functions and so on). In contrast, the study of complexity suggests to consider first more restricted classes of objects admitting fast enough computations. Here we briefly (and not very systematically) recall some relevant algebraic notions and facts. For details see e.g. \cite{wa67}.

Coefficients of matrices and polynomials will usually be taken from a fixed (ordered) field $\mathbb{F}\subseteq\mathbb{R}$ of reals. For a field $\mathbb{F}$, let $\mathbb{F}[x_1,x_2,\ldots]$ denote the ring of polynomials with coefficients in $\mathbb{F}$ and variables $x_1,x_2,\ldots$, and let $\mathbb{F}(x_1,x_2,\ldots)$ be the corresponding field of fractions (i.e., the field of rational functions with variables $x_1,x_2,\ldots$). A polynomial $p\in\mathbb{F}[x]$ is written as $p=a_0+a_1x+\cdots+a_nx^n$ where $a_n\not=0$ for $p\not=0$; $n$ is the {\em degree} of $p$ denoted as ${\rm deg}(p)$; if $a_n=1$, $p$ is called {\em unitary}. A  polynomial $p\in\mathbb{F}[x]$ is {\em irreducible} (over $\mathbb{F}$) if it is not the product of two polynomials from $\mathbb{F}[x]$  of lesser degrees. 

The arithmetics of polynomials is very similar to the arithmetics of the integers $\mathbb{Z}$, in particular any non-zero polynomial $p$, ${\rm deg}(p)\geq2$, has a canonical (i.e., unique up to permutation of factors) factorisation $p=ap_1^{m_1}\cdots p_k^{m_k}$ where $a\in\mathbb{F}$, $k,m_1,\ldots,m_k\geq1$ and $p_1,\ldots,p_k$, ${\rm deg}(p_i)\geq1$, are unitary irreducible  polynomials such that $(p_i,p_j)=1$ for $i\not=j$ ($(p_i,p_j)$ is the greatest common divisor of $p_i,p_j$). It is known that if $0\not=\alpha\in\mathbb{C}$ is algebraic over $\mathbb{F}$ (i.e. $\alpha$ is a root of some polynomial over $\mathbb{F}$) then the smallest subfield $\mathbb{F}(\alpha)$ of $\mathbb{C}$ containing $\mathbb{F}\cup\{\alpha\}$ is isomorphic to the quotient field $\mathbb{F}[x]/(p_\alpha)$ where $p_\alpha$ is the unique unitary polynomial of minimal degree $p\in\mathbb{F}[x]$  with $p(\alpha)=0$ and $(p_\alpha)$ is the ideal of $\mathbb{F}[x]$ generated by $p_\alpha$; ${\rm deg}(\alpha)={\rm deg}(p_\alpha)$ is called the {\em degree of $\alpha$ over} $\mathbb{F}$.

In this paper we most often work  with  ordered fields $\mathbb{F}\in\{\mathbb{Q},\mathbb{Q}(\alpha),\mathbb{A}\mid 0\not=\alpha\in\mathbb{A}\}$ where $\mathbb{Q}$ is the ordered field of rationals, $\mathbb{A}$ is the ordered field of algebraic reals (which consists of the reals algebraic  over $\mathbb{Q}$).  One can also consider the smallest subfield $\mathbb{Q}(\alpha_1,\ldots,\alpha_n)$ of $\mathbb{R}$ containing given $\alpha_1,\ldots,\alpha_n\in\mathbb{A}$. For any such field there is a ``primitive element'' $\alpha\in\mathbb{A}$ with $\mathbb{Q}(\alpha)=\mathbb{Q}(\alpha_1,\ldots,\alpha_n)$. Note that $\mathbb{A}$ is the smallest real closed ordered field. With any non-zero $\alpha\in\mathbb{A}$ we associate the unique pair $(p_\alpha,k)$ such that $k$ satisfies $\alpha=\alpha_k$ where $\alpha_1<\cdots<\alpha_m$ is the increasing sequence of all real roots of the  polynomial $p_\alpha$.

Let $M_n(R)$ be the set of $n\times n$-matrices over a (commutative associative with a unit element 1) ring $R$, and $M(R)$ be the union of all $M_n(R)$, $n\geq1$. We use standard terminology and notation from linear algebra. In particular, ${\rm det}(A)$ is the determinant of $A=(a_{ij})\in M_n(R)$, ${\rm diag}(a_1,\ldots,a_n)$ is the diagonal matrix with the diagonal elements $a_1,\ldots,a_n\in R$, so in particular $I=I_n={\rm diag}(1,\ldots,1)$ is the unit matrix. The roots of the polynomial $ch_A={\rm det}(\lambda I-A)$ are called {\em eigenvalues of} $A\in M_n(\mathbb{C})$. In general, the eigenvalues of a real matrix are complex numbers. The eigenvalues of a symmetric real matrix are always real.

\subsection{Encodings and bit complexity}\label{coding}

Computations on  existing computers (as well as on theoretical computing devices like Turing machines) don't work  with abstract mathematical objects (like integers, rationals or polynomials) but rather with words over a finite alphabet. Here we briefly recall some relevant notions and facts (for more details see e.g. \cite{bdg,sc86,lo82,cr91,al16,as17}).

Let $\Sigma$ be a finite alphabet, $\Sigma^\ast$  the set of words over $\Sigma$, $S\subseteq(\Sigma^\ast)^n$, and $t:S\to\omega$. A  function $f:S:\to\Sigma^\ast$ is {\em computable in time $t$} if there exists a $k$-tape Turing machine, $k\geq n+1$, such that starting to work with words $x_1,\ldots,x_n=\bar{x}\in S$ written on the first $n$ tapes, finishes within $t(\bar{x})$ steps with $f(\bar{x})$ written on the $(n+1)$-st tape. The set $S$ is {\em computable in time $t$} if so is its characteristic function $\chi_S:(\Sigma^\ast)^n\to\{0,1\}$ (assuming that $0,1\in\Sigma$). A structure $\mathbb{S}=(S;\ldots)$ of a finite signature is  computable in time $t$ if so are its universe $S$ and all the signature functions and relations.  An abstract structure $\mathbb{A}$ is {\em $t$-time-presentable} if it is isomorphic to a structure $\mathbb{B}$ computable in time $t$. Any isomorphism from $\mathbb{A}$ onto $\mathbb{B}$ is a {\em $t$-time-presentation} of $\mathbb{A}$. 
Note that usually people work not directly with  Turing machines but rather with informal algorithms on words; the algorithms should use only elementary enough steps to make it clear how to translate them to the syntax of Turing machines.

If the function $t(\bar{x})$  is bounded by some polynomial on $|\bar{x}|={\rm max}\{|x_1|,\ldots,|x_k|\}$ (where $|x|$ is the length of a word $x$) we say that the function $f$ (resp., the set $S$) is computable in polynomial time ({\em $p$-computable} for short). An abstract structure $\mathbb{A}$ is {\em $p$-presentable} if it is isomorphic to a $p$-computable structure $\mathbb{B}$. Any isomorphism from $\mathbb{A}$ onto $\mathbb{B}$ will be called a {\em $p$-presentation} of $\mathbb{A}$. Similarly one can define computability in linear, quadratic or  exponential time.  

In order to enable Turing machines  work with abstract objects, we have to encode the objects by words in a finite alphabet. First note (see e.g. \cite{bdg}) that for  any finite alphabet $\Sigma$ there is a natural encoding (injective function) $c:\Sigma^\ast\to\{0,1\}^\ast$ of words over $\Sigma$ by binary words such that the function $c$, its range $rng(c)$, and the inverse function $c^{-1}$ are computable in linear time, so in all interesting cases we can without loss of generality stick to  binary encodings (this is the reason why the classical computational complexity is often called bit complexity). We give several examples of such encodings for different sets $A$ of abstract objects. In describing such an encoding $b:A\to\{0,1\}^\ast$ we often use the  trick of first describing an auxiliary encoding $e:A\to\Sigma^\ast$ for some bigger alphabet $\Sigma$ and then setting $b=c\circ e$. Below we often take $\Sigma=\{0,1,\ast\}$ where $\ast$ is a new symbol used for separating binary words.

A positive integer $n$ is usually encoded by its binary notation $b(n)$, so $|b(n)|=log (n)$. Adding an additional bit for the sign, we obtain a binary encoding $b$ of the integers. Identifying rationals with fractions $\frac{p}{q}$ where $p,q\in\mathbb{Z} $, $q\geq1$ and $(p,q)=1$, we can define an encoding $c:\mathbb{Q}\to\{0,1,\ast\}^\ast$ by $c(\frac{p}{q})=b(p)\ast b(q)$ (and, by the mentioned trick with alphabets, we can modify $c$ to obtain a binary coding $b$ of the rationals). The defined encodings give $p$-presentations of the ordered ring $\mathbb{Z}$ and the ordered field $\mathbb{Q}$.

The binary encoding $b$ of $\mathbb{Q}$ induces the encoding  $e:\mathbb{Q}[x]\to\{0,1,\ast\}^\ast$ which associates with a non-zero polynomial $p=a_0+a_1x+\cdots+a_nx^n$, $a_n\not=0$, the code $b(a_0)\ast\cdots\ast b(a_n)$. In a similar way one can define natural induced encodings of $\mathbb{Q}[x_1,\ldots,x_n]$ and of $\mathbb{Q}(x_1,\ldots,k_n)$ for each $n\geq1$ which provide $p$-presentations of the corresponding rings and fields (see e.g. \cite{as17} for additionl details). Moreover, in the field $\mathbb{Q}(x_1,\ldots,x_n)$ also the subtraction and division, as well as the evaluation  function $\mathbb{Q}[x_1,\dots,x_n]\times\mathbb{Q}^n\to\mathbb{Q}$ are $p$-computable. 
Furthermore, these encodings induce  encodings of $\mathbb{Q}[x_1,\dots]=\bigcup_n\mathbb{Q}[x_1,\dots,x_n]$ and of $\mathbb{Q}(x_1,\dots)$ modulo which any of $\mathbb{Q}[x_1,\dots,x_n],\mathbb{Q}(x_1,\dots,x_n)$ is $p$-computable; moreover,
the evaluation partial function $\mathbb{Q}[x_1,\ldots]\times\mathbb{Q}^\ast\to\mathbb{Q}$ (and the similar function for $\mathbb{Q}(x_1,\ldots)$), where $\mathbb{Q}^\ast$ is the set of finite strings over $\mathbb{Q}$, is $p$-computable.

If $p\in\mathbb{Q}[x]$ is  unitary irreducible then the quotient-field $\mathbb{Q}[x]/(p)$ (formed by the polynomials of degree $<{\rm deg}(p)$)  also has a natural $p$-presentation (see e.g. \cite{ak89,as17}). Therefore, the ordered field $\mathbb{Q}(\alpha)$, for each non-zero $\alpha\in\mathbb{A}$, has a natural $p$-presentation.
For the induced binary presentation of $\mathbb{F}=\mathbb{Q}(\alpha)$ we can repeat the constructions of the previous paragraph and obtain natural $p$-presentations of $\mathbb{F}[x_1,\dots]$ and of $\mathbb{F}(x_1,\ldots)$. Again, the polynomial evaluation $\mathbb{F}[x_1,\dots]\times\mathbb{F}^\ast\to\mathbb{F}$ will be $p$-computable w.r.t. the natural induced encodings. 

We also define a natural encoding of $\mathbb{A}$ into $\{0,1,\ast\}^\ast$ by associating with any non-zero $\alpha\in\mathbb{A}$ the word $b(p_\alpha)\ast b(k)$ where $(p_\alpha,k)$ is the pair from the previous subsection and $b:\mathbb{Q}[x]\to\{0,1\}^\ast$ is the natural binary coding specified above. From deep results of computer algebra (including the polynomial time algorithm for factoring rational polynomials \cite{lll82}) it follows that in this way we obtain a natural $p$-presentation of the ordered field $\mathbb{A}$ in which also subtraction and division are $p$-computable (see e.g. \cite{lo82,as17}). 

Simplifying notation, we always denote the natural binary $p$-presentations of any of the ordered fields $\mathbb{F}\in\{\mathbb{Q},\mathbb{Q}(\alpha),\mathbb{A}\mid 0\not=\alpha\in\mathbb{A}\}$   by $b$. Note however that some important computational properties of the presentation of $\mathbb{A}$ discussed above differ from those for the presentations of $\mathbb{Q}$ and $\mathbb{Q}(\alpha)$. In particular, the evaluation function $\mathbb{A}[x_1,\dots]\times\mathbb{A}^\ast\to\mathbb{A}$ is now computable in exponential time but  not in polynomial time. In fact, already the  ``long sum'' operation $(\alpha_1\ast\cdots\ast\alpha_n)\mapsto\alpha_1+\cdots+\alpha_n$ is not PTIME-computable uniformly on $n$ (even not computable in PSPACE) w.r.t. the presentation of $\mathbb{A}$. This follows from the results  in \cite{zhou} (for a detailed explanations see comments after the proof of Theorem 4 in \cite{as17}).

We will use some results from \cite{as17} about the complexity of root-finding in the field $\ca=(\caa;+,\times,0,1)$ of complex algebraic numbers, i.e. of finding all roots of an equation
$ \alpha_{e}x^{e} + \ldots + \alpha_{1}x + \alpha_{0} = 0 $ where
$ \alpha_{i} \in \caa $ for $ i \ls e $. 
More precisely, the authors of  \cite{as17} consider equations of the form
\begin{equation}\label{eq}
 t_{e} (\alpha_{1}, \ldots, \alpha_{k})x^{e} + \ldots +
t_{1} (\alpha_{1}, \ldots, \alpha_{k})x + t_{0} (\alpha_{1}, \ldots, \alpha_{k}) = 0,
 \end{equation}
 where $ \alpha_{1}, \ldots, \alpha_{k} \in \caa $ and
$ t_{j} (\bx) \in \rat [x_{1}, \ldots, x_{k}] $. The problem is to find  a list of (codes of) all roots of (\ref{eq}) from given $b( \alpha_{1}) * \cdots * b(\alpha_{k}) $ and $ b(t_{0} (\bx)) * \cdots * b(t_{e} (\bx))$ where $b$ is a natural binary encoding of $\ca$ induced by the presentation of $\mathbb{A}$ described above and by Gauss representation of complex numbers as pairs of reals. As shown in Theorem 8  \cite{as17}, the problem is solvable in polynomial time  for any fixed $k$. Moreover, the same estimate holds for the version of this problem when one computes the list of all distinct real roots of (\ref{eq}) in increasing order.

The introduced binary presentations of fields, polynomials and rational functions are natural in the sense that they are $p$-equivalent to some presentations really used in computer algebra systems. See e.g. \cite{lo82,as17} for additional details.

 Associate with any matrix $A=(a_{ij})\in M_n(\mathbb{F})$ its code $c(A)\in\{0,1,\ast\}^\ast$ by $c(A)=b(a_{11})\ast\cdots \ast b(a_{1n})\ast b(a_{21})\ast\cdots\ast b(a_{2n})\ast\cdots\ast b(a_{n1})\ast\cdots\ast b(a_{nn})$. These encodings induce the binary encoding of the set $M(\mathbb{F})=\bigcup_nM_n(\mathbb{F})$ of all square matrices over $\mathbb{F}$ in which any set $M_n(\mathbb{F})$ is $p$-computable. Furthermore, many matrix properties like symmetricity are also $p$-computable. 
 
It is well known (see e.g. \cite{sc86}) that our encodings give $p$-presentations of the rings $M_n(\mathbb{Q})$ uniformly on $n$ (uniformity means that there is a polynomial bound working for all $n$). Moreover, it is easy to check that evaluation of some ``long'' terms in these rings are also $p$-computable uniformly on $n$ w.r.t. these presentations (in particular the function $A_1\ast\cdots\ast A_n\mapsto A_1\times\cdots\times A_n$ is $p$-computable). Even more involved matrix algorithms like computing of the determinant, computing of the inverse of a non-degenerate matrix, and Gauss method also work in polynomial time (see Chapter 3 \cite{sc86} for additional details). 

From results in \cite{lo82,as17} (using also the arguments in \cite{sc86}) it follows that the polynomial time estimates of the previous paragraph remain true for the rings $M_n(\mathbb{Q(\alpha)})$ for each non-zero $\alpha\in\mathbb{A}$. Moreover, our presentation of the ring $M_n(\mathbb{A})$ is a $p$-presentation for any fixed $n\geq1$, but not uniformly on $n$.  The ``long terms'' are not $p$-computable w.r.t. our presentation for  $M_n(\mathbb{A})$, even for a fixed $n$.

The introduced encodings of matrices are natural in the sense they are closely related to standard encodings of matrices in numerical  analysis. The only difference is that here we use a precise symbolic encoding of matrix coefficients while in numerical analysis the floating-point approximations of coefficients are usually used.

\begin{Remark}\label{complex_remark}
Note that, along with  bit complexity, there are other approaches to measuring  complexity of computations on, say, Euclidean spaces, including algebraic complexity (see e.g. \cite{ber}),  topological complexity (see e.g.  \cite{sma,vas}), real complexity, i.e. bit complexity adapted to reals by means of computable analysis  (see e.g. \cite{kz,bhw,wei03,wei}). Algebraic complexity counts the arithmetical operations needed for finding a symbolic solution (the notion we use is close to that, but also taking into account the bit costs of performing the operations). Topological complexity counts the number of equality and order tests needed to find a solution of a discontinuous (hence, non-computable in the sense of \cite{bhw,wei}) problem. Real complexity  counts steps of a Turing machine working on the representations of reals (when restricted to algebraic reals, it is close to our approach). \end{Remark}

\subsection{Formulations of main results}\label{main_subsec}

Now we have enough notions and terminology to state the guaranteed-precision problems in a rigorous form.

First we consider  the task of computing the domain $H$ of existence and uniqueness of the Cauchy problem. As mentioned in Section \ref{cbv}, the set $H$ is the intersection of  semi-spaces
 $$
t\geq0,\;x_i-\mu^{(i)}_{\rm
max}t\geq0,\;x_i-1-\mu^{(i)}_{\rm min}t\leq0,\;(i=1,\ldots,m)
 $$
of ${\mathbb R}^{m+1}$ where $\mu^{(i)}_{\rm
min},\mu^{(i)}_{\rm
max}$ are respectively the minimum and maximum of the eigenvalues of the matrix  $A^{-1}B_i$. Therefore, the computation of $H$ reduces to the computation of the eigenvalues of the matrices  $A^{-1}B_i$. 

Our algorithms for solving the Cauchy problem will be for technical reasons presented only for the case when $H$ satisfies the condition $\mu^{(i)}_{\rm
min}<0<\mu^{(i)}_{\rm max}$ for all $i=1,\ldots,m$; this condition
often holds for natural physical systems.  Note that this condition implies that $H$ is a
compact subset of $[0,+\infty)\times Q$.

In \cite{ss09} we observed that the domain $H$ for the problem \eqref{sist_1}
is computable from matrices $A,B_1,\ldots,B_m$ (more exactly, the vector
$(\mu^{(1)}_{\max},\ldots,\mu^{(m)}_{\max},
\mu^{(1)}_{\min},\ldots,\mu^{(m)}_{\min})$  is computable
from $A,B_1,\ldots,B_m$; this implies the computability of $H$ in
the sense of computable analysis \cite{wei}).

The  next result establishes the complexity of computing $H$ satisfying the mentioned condition. 

\begin{Theorem}\label{cpH}
Let  $m,n\geq2$ be any fixed integers. There is a polynomial time algorithm which for any given $A,B_1\ldots,B_m\in M_n(\mathbb{A})$ finds the vector
$(\mu^{(1)}_{\max},\ldots,\mu^{(m)}_{\max},
\mu^{(1)}_{\min},\ldots,\mu^{(m)}_{\min})$ and checks the condition $\mu^{(i)}_{\rm
min}<0<\mu^{(i)}_{\rm max}$ for all $i=1,\ldots,m$. Thus, the algorithm finds the domain $H$ satisfying the condition above, or reports on the absence of such a domain.
 \end{Theorem}

Now we state two  guaranteed-precision versions of a restricted (to algebraic matrix coefficients and rational polynomials as initial data) Cauchy problem  \eqref{sist_1}. Let $m,n$  be  fixed positive integers. 

{\bf Task 1}. CP$(m,n,\mathbb{A},\mathbb{Q})$ is the following computational task:

INPUT: Integer $a\geq1$,
polynomials $\varphi_1\ldots,\varphi_n\in \mathbb{Q}[x_1\ldots,x_m], \;f_1\ldots,f_n\in \mathbb{Q}[t,x_1\ldots,x_m]$ and matrices $A,B_1\ldots,B_m\in M_n(\mathbb{A})$.

OUTPUT: (Codes of) a rational $T>0$ with $H\subseteq[0,T]\times Q$, a spatial rational grid step $h$ dividing $1$, a time grid step $\tau$ dividing $T$ and a rational $h,\tau$-grid function $v:G_N^{\tau}\to\mathbb{Q}$ 
such that \begin{equation}\label{guarant}||{\bf u}-\widetilde{\upsilon\mid_H}||_{sL_2}<\frac{1}{a}.\end{equation} 
  
This problem is a standard computational task on a discrete structure which  asks  for an algorithm (and its complexity estimation) which, for any given input computes a suitable output. 
The second version of the Cauchy problem is formulated as follows.

{\bf Task 2.} CP$(m,n,a,M,\mathbb{A},\mathbb{Q})$ is the following computational task:

INPUT: (Codes of) polynomials $$\varphi_1\ldots,\varphi_n\in \mathbb{Q}[x_1\ldots,x_m], \;f_1\ldots,f_n\in \mathbb{Q}[t,x_1\ldots,x_m]$$ and matrices $A,B_1\ldots,B_m\in M_n(\mathbb{A})$ such that the  quantities $||A||_2$, $\frac{\lambda_{max}(A)}{\lambda_{min}(A)}$,  $$\max_i\Bigl\{||B_i||_2,  
 ||(A^{-1}B_i)^2||_2,  \max_k\{|\mu_k|:\text{det}(\mu_k A-B_i)=0\}, \sup_{t,x}||\frac{\partial^2
f}{\partial x_i\partial t}(t,x)||_2 \Bigr\}, $$ and
  $$\max_{i,j}\Bigl\{||A^{-1}B_iA^{-1}B_j-A^{-1}B_jA^{-1}B_i||_2,
\sup_{t,x}||\frac{\partial^2
f}{\partial x_i\partial x_j}(t,x)||_2,  \sup_x||\frac{\partial^2
\varphi}{\partial x_i\partial x_j}(x)||_2\Bigr\}$$ are bounded by $M$.

OUTPUT: (Codes of)a rational $T>0$ with $H\subseteq[0,T]\times Q$, a spatial rational grid step $h$ dividing $1$, a time grid step $\tau$ dividing $T$ and a rational $h,\tau$-grid function $v:G_N^{\tau}\to\mathbb{Q}$ 
such that
 $||{\bf u}-\widetilde{\upsilon\mid_H}||_{sL_2}<\frac{1}{a}.$

The guaranteed-precision versions BV$(m,n,a,M,\mathbb{A},\mathbb{Q})$ and BV$(m,n,\mathbb{A},\mathbb{Q})$ of  boundary-value problem  \eqref{sist_2} are stated in a similar way, only in the second version  the additional assumption $\max_i\{||\Phi_i^{(1)}||_2,\ ||\Phi_i^{(2)}||_2\}\leq M$ is needed. 

Our basic result, which concerns Tasks 1 and 2, is formulated as follows (the definitions of the well known complexity classes PTIME and EXPTIME may be found e.g. in \cite{bdg}).

\begin{Theorem}\label{cpq}
\begin{enumerate}
 \item  For any  $m,n\geq 1$, the problems CP$(m,n,\mathbb{A},\mathbb{Q})$  and  BVP$(m,n,\mathbb{A},\mathbb{Q})$  are solvable in EXPTIME. 
 \item  For any  $m,n,a,M\geq1$, the problems CP$(m,n,a,M,\mathbb{A},\mathbb{Q})$ and   BVP$(m,n,a,M,
\mathbb{A},\mathbb{Q})$  are solvable in PTIME. 
%\item  For any  $m,a,M\geq1$, the problems CP$(m,a,M,\mathbb{Q},\mathbb{Q})$ and  BVP$(m,a,M,\mathbb{Q},\mathbb{Q})$ from Task 2 are solvable in PTIME. 
\end{enumerate}
\end{Theorem}

One could consider several  variations of the stated problems.
For instance, one could vary some of the remaining fixed parameters or take rational functions $\varphi_1\ldots,\varphi_n\in \mathbb{Q}(x_1\ldots,x_m), \;f_1\ldots,f_n\in \mathbb{Q}(t,x_1\ldots,x_m)$ instead of polynomials. We will make some comments on the first option below. For the second option, analogues of Theorem \ref{cpq} may be proved using the methods of this paper. In our approach, we cannot take arbitrary algebraic polynomials as the initial and right-hand side functions because computing their values, by the above-mentioned result in \cite{zhou}, requires exponential time which yields worse upper bounds than in Theorem \ref{cpq} (see Subsection \ref{csteps}).

\begin{Remark}
From the proof in the last section of this paper it will follow that similar upper bounds hold for the algebraic complexity.
\end{Remark}

\section{Auxiliary algorithms}\label{aux}

In this section we present some linear algebra algorithms which are  used in Section \ref{mainres} to prove  the main results.

\subsection{Computing spectral decomposition}\label{spectral}

Here we give upper bounds for the complexity of symbolic computations of  eigenvalues and eigenvectors for some classes of matrices and matrix pencils. The computations are w.r.t. the encodings specified in Section \ref{coding}.

By {\em spectral decomposition} of a symmetric real matrix $A\in M_n(\mathbb{R})$ we mean a pair $((\lambda_1,\ldots,\lambda_n),(\mathbf{v}_1,\ldots,\mathbf{v}_n))$ where $\lambda_1\leq\cdots\leq\lambda_n$ is the non-decreasing sequence of all eigenvalues of $A$ (each eigenvalue occurs in the sequence several times, according to its multiplicity) and $\mathbf{v}_1,\ldots,\mathbf{v}_n$ is a corresponding orthonormal basis of eigenvectors. 

\begin{Proposition}\label{specmat}
\begin{enumerate}
 \item For any fixed $n\geq1$, there is a polynomial time algorithm which, given a symmetric matrix $A\in M_n(\mathbb{A})$, computes a spectral decomposition of $A$.
 
 \item There is an algorithm which, given a symmetric matrix $A\in M(\mathbb{Q})$, computes (uniformly on $n$) the spectrum $(\lambda_1,\ldots,\lambda_n)$ of $A$ in PTIME and a spectral decomposition of $A$ in EXPTIME.\footnote{Formulation of this item is weaker than in \cite{ss18} where the estimate for orthonormalized eigenvectors is incorrect because the Gram-Schmidt process does not work in PTIME uniformly on $n$.}
\end{enumerate}
\end{Proposition}

{\em Proof.} 1. For fixed $n$, the coefficients of $ch_A={\rm det}(\lambda I-A)\in\mathbb{A}[\lambda]$ are fixed integer polynomials from the elements of $A$, so we can compute them in polynomial time using Theorem 4 in \cite{as17}. 
By results in \cite{lo82} and Theorem 8 in \cite{as17} cited in Section \ref{algprel}, we can compute in polynomial time the increasing sequence $\mu_1<\cdots<\mu_m$ of all roots of $ch_A$ and the corresponding multiplicities $r_1,\ldots,r_m$ (hence we can also compute the sequence $\lambda_1,\ldots,\lambda_n$). Compute in polynomial time a primitive element $\beta\in\mathbb{A}$ with $\mathbb{Q}(\beta)=\mathbb{Q}(\alpha,\mu_1,\ldots,\mu_m)$. Working in $\mathbb{Q}(\beta)$, we find in polynomial time (in the usual way, solving corresponding linear systems by the Gauss method), for each $j=1,\ldots,m$,  a basis $\mathbf{w}^i_1,\ldots,\mathbf{w}^i_{r_j}$ for the eigenspace of $\mu_j$. Applying the Gram-Schmidt orthogonalisation process (which also works in polynomial time for fixed $n$) and normalizing the obtained orthogonal basis, we obtain an orthonormal basis for this eigenspace. Putting together the orthonormal bases for all $j$, we obtain a desired orthornormal basis $(\mathbf{v}_1,\ldots,\mathbf{v}_n)$ for the whole space.  

2. Given $A\in M(\mathbb{Q})$, compute first the order $n$ of $A$.  Next we compute the characteristic polynomial $ch_A={\rm det}(\lambda I-A)\in\mathbb{Q}[\lambda]$, written in the form $ch_A=\lambda^n-p_1\lambda^{n-1}-p_2\lambda^{n-2}-\cdots-p_n$. By remarks at the end of Section \ref{coding}, we can compute in polynomial time  the traces $s_1,\ldots,s_n$ of matrices $A^1,\ldots,A^n$ respectively. It is known (see e.g. Section 4.4 of \cite{gant}) that the equalities $p_1=s_1$, and $(k+1)p_{k+1}=s_{k+1}-p_ks_k$ for $k=1,\ldots,n-1$, hold. From these we subsequently compute $p_1,p_2,\ldots,p_n$ in polynomial time. 

Since factorisation in $\mathbb{Q}[\lambda]$ is computable in polynomial time \cite{lll82}, we can compute the canonical decomposition $ch_A= q_1^{k_1}\cdots q_l^{k_l}$ where $l,k_1,\ldots,k_l$ are positive integers. Let $\lambda_{j,1}<\cdots<\lambda_{j,d_j}$ be the sequence of all roots of $q_j$ where $d_j={\rm deg}(q_j)$. Then $c(q_j,i)$ is the code of $\lambda_{j,i}$ for each $i=1,\ldots,d_j$. Since $\{\lambda_{j,i}\mid j=1,\dots,l,\;i=1,\ldots,d_j\}=\{\lambda_1,\ldots,\lambda_n\}$ and all $\lambda_{j,i}$ are pairwise distinct, we have computed the (codes of the) eigenvalues $\lambda_1,\ldots,\lambda_n$. Since the multiplicity of $\lambda_{j,i}$ is $k_j$, we have also computed the multiplicity $r_i$ of any eigenvalue $\lambda_i$ (the multiplicity of $\lambda_{j,i}$ is $k_j$).

Since the eigenvalues $\lambda_{j,i}$ are pairwise distinct, it remains to find an orthonormal basis of the eigenspace $\{\mathbf{x}\mid (\lambda^i_j\cdot I-A)\mathbf{x}=\mathbf{0})\}$ corresponding to any fixed $\lambda_{j,i}$ (the dimension of this space is $k_j$). Applying the Gauss method (in the field of coefficients $\mathbb{Q}(\lambda_{j,i}$) we find a desired basis $\mathbf{w}_{j,1},\ldots,\mathbf{w}_{j,k_j}$ (note that, since the field $\mathbb{Q}(\lambda_{j,i})$ is isomorphic to $\mathbb{Q}[\lambda]/(p_j)$ for each $i=1,\ldots,d_j$, the systems may be solved uniformly on $i$ in the quotient field, with subsequent evaluation of the computed polynomials at $\lambda_{j,i}$). By remarks in Section \ref{coding}, this computation runs in polynomial time. Applying the Gram-Schmidt orthogonalization process and normalizing the obtained orthogonal bases for each $i=1,\ldots,d_j$, we put them together to get a resulted orthonormal basis for the whole space in EXPTIME uniformly on $n$, using Theorems 4 and 8 in \cite{as17}. 
 $\qed$
 
 \begin{Remark}\label{rem_spec}
 In fact, Theorem 8 in \cite{as17} implies that a spectral decomposition of any given symmetric algebraic matrix may be computed in exponential time uniformly on $n$, though in this case we cannot compute in polynomial time even coefficients of $ch_A$, again because of the result in \cite{zhou}.
  \end{Remark}

By  {\em matrix pencil} we mean a pair $(A,B)$ (often written in the form $\mu A-B$) of real non-degenerate symmetric matrices such that $A$ is positive definite (i.e., all of its eigenvalues are positive).
By {\em spectral decomposition} of such a pencil we mean a tuple
 $$((\lambda_1,\ldots,\lambda_n),(\mathbf{v}_1,\ldots,\mathbf{v}_n),(\mu_1,\ldots,\mu_n),(\mathbf{w}_1,\ldots,\mathbf{w}_n))$$
  such that $((\lambda_1,\ldots,\lambda_n),(\mathbf{v}_1,\ldots,\mathbf{v}_n))$ and $((\mu_1,\ldots,\mu_n),(\mathbf{w}_1,\ldots,\mathbf{w}_n))$ are spectral decompositions  of the symmetric matrices $A$  and $D^*L^*BLD$ respectively, where  $L$ is the matrix formed by vectors $\mathbf{v}_1,\ldots,\mathbf{v}_n$ written as columns and $D=\operatorname{diag}\{\frac{1}{\sqrt{\lambda_1}}, \frac{1}{\sqrt{\lambda_2}}, \ldots, \frac{1}{\sqrt{\lambda_n}}\}$.
  
\begin{Proposition}\label{specpen}
For any fixed $n\geq1$, there is a polynomial time algorithm which, given a matrix pencil $(A,B)$ with $A,B\in M_n(\mathbb{A})$, computes a spectral decomposition of $(A,B)$.
\end{Proposition}

{\em Proof.}  By item 1 of the previous theorem, we can find in polynomial time a spectral decomposition $((\lambda_1,\ldots,\lambda_n),(\mathbf{v}_1,\ldots,\mathbf{v}_n))$ of $A$. Since we can solve polynomial equations in $\mathbb{Q}(\lambda_1,\ldots,\lambda_n)$ in polynomial time (see e.g. \cite{lo82,as17} for details), we can compute in polynomial time the matrix $D^*L^*BLD$. Applying item 1 of the previous theorem to this matrix (and working now in the field $\mathbb{Q}(\lambda_1,\ldots,\lambda_n)$ which is also computable in polynomial time \cite{lo82}), we compute the remaining items $(\mu_1,\ldots,\mu_n),(\mathbf{w}_1,\ldots,\mathbf{w}_n)$. Note that $(\mu_1,\ldots,\mu_n)$ coincides with the sequence of eigenvalues of (in general, non-symmetric) matrix $A^{-1}B$. 
 $\qed$

\subsection{Computing data for the difference scheme}

Here we explain how to compute  data needed for  computations with the difference schemes in Section \ref{mainres}. 

Let $A,B_1,\dots,B_m\in M_n(\mathbb{A})$ be matrices satisfying the conditions in Cauchy problem. We can compute the spectral decomposition $((\lambda_1,\ldots,\lambda_n),(\mathbf{v}_1,\ldots,\mathbf{v}_n))$  of $A$  as in the proof of Proposition \ref{specmat}. Let $\lambda_{max}$, $\lambda_{min}$ be respectively the maximum and minimum of $\lambda_1,\ldots,\lambda_n$. Let $L$ be the orthonormal matrix formed by vectors $\mathbf{v}_1,\ldots,\mathbf{v}_n$ written in columns, so
 $L^*AL=\Lambda={\rm diag}\{\lambda_1,\lambda_2,\ldots,\lambda_n\}$, and let
$D=\Lambda^{-\frac{1}{2}}$. 

For each $i=1,\ldots,m$, let $((\mu^{(i)}_1,\ldots,\mu^{(i)}_n),(\mathbf{w}^{i}_1,\ldots,\mathbf{w}^{i}_n))$ be the spectral decomposition  of the symmetric matrix $D^*L^*B_iLD$ computed as in the proof of Proposition \ref{specpen}. Let $\mu^{(i)}_{max}$, $\mu^{(i)}_{min}$ be respectively the maximum and minimum of $\mu^{(i)}_1,\ldots,\mu^{(i)}_n$. Let $M_i={\rm
 diag}\{\mu^{(i)}_1,\ldots,\mu^{(i)}_n)\}$ and  $K_i$ be the orthonormal matrix formed by vectors $\mathbf{w}^{i}_1,\ldots,\mathbf{w}^{i}_n$ written in columns, so $K_i^*D^*L^*B_iLDK_i=M_i$. Let $T_i=LDK_i$ for each $i=1,\ldots,m$.

From Propositions \ref{specmat}, \ref{specpen} and remarks in Section \ref{coding} we easily obtain:

\begin{Proposition}\label{schemedata}
For any fixed $m,n\geq 1$, there is a polynomial time algorithm which, given  matrices $A,B_1,\dots,B_m,\in M_n(\mathbb{Q})$ satisfying the conditions of symmetric hyperbolic systems, computes the objects $A^{-1}$, $T_i$, $T_i^{-1},\lambda_{max}$, $\lambda_{min},\mu^{(i)}_{max},\mu^{(i)}_{min}$, $\mu^{(i)}_k$($i=1,\ldots,m,k=1,\ldots, n$) specified above. 
\end{Proposition}

Note that in the proof of Theorem \ref{cpq} in Section \ref{mainres} we will stick (for notational simplicity) to the typical particular case  $m=2$ where the notations $T_x=T_1,T_y=T_2$ and $K_x=K_1,K_y=K_2$ are more appealing, e.g. for considering the linear transformations of variables.

\section{Proof of the main results} \label{mainres}

Theorem \ref{cpH} straightforwardly follows from Proposition \ref{specpen}. We give a proof of Theorem \ref{cpq}  for the Cauchy problem and the boundary-value problem simultaneously, with a minor modification of the numerical algorithm used in our proof.

\subsection{Computing the grid steps}\label{steps}

In the Proposition below we assume that all the assumptions of Section \ref{main_subsec}  hold. This proposition is the main technical tool of the present paper.

\begin{Proposition}\label{steps_prop} Let $m,n\geq 1$ be fixed integers. The time and space steps $\tau$ and $h$  guaranteeing \eqref{guarant}, when using the difference scheme \eqref{scheme_2_sist1} on the corresponding grid, are  PTIME computable from $a, \varphi, f, A, B_i$ ($i=1,2,\ldots,m$). \end{Proposition}

{\em Proof.} W.l.o.g., stick to the case $m=2$, denoting $B_1=B$, $B_2=C$. First note \cite{go76} that the Godunov scheme (see its detailed description below in the next subsection)  is stable  if and only if 
\begin{equation}\label{tau}\tau\leq h\cdot \left(\frac{1}{\max_i\{|\mu_i|:\text{det}(\mu_i A-B)=0\}}+\frac{1}{\max_i\{|\mu_i|:\text{det}(\mu_i A-C)=0\}}\right)^{-1}\end{equation}
(see also a short summary of the proof of this fact in \cite{ss09}). Recall that stability is an intrinsic property of a difference scheme, implying, together with approximation, its convergence to the corresponding differential equation in grid norms.

 By Proposition \ref{specpen}, a rational $\tau$ satisfying \eqref{tau} can be found in polynomial time, if $h$ is computed in polynomial time. It is also obvious that $\tau$ can be chosen so that $L=\frac{T}{\tau}$ is integer.  Thus it suffices to estimate complexity of finding $h$.

In \cite{ss13} we established, applying well known theorems about difference schemes and interpolations, that 
 \begin{equation}\label{est1} 
 ||{\bf u}-\widetilde{\upsilon\mid_H}||_{sL_2}\leq ||{\bf u}-\widetilde{{\bf u}\mid_{G_N^{\tau}}}||_{sL_2}+||\widetilde{{\bf u}\mid_{G_N^{\tau}}}-\widetilde{\upsilon\mid_H}||_{sL_2}
\leq 
c_{int}h+c_{diff}h\leq {\cal C}h\leq \frac{1}{a},\end{equation}
 where the ``constants'' $c_{int}$ and $c_{diff}$ depend only on $A, B_1,\ldots, B_m$ and the first and second partial derivatives of $\varphi$.
So we can take $h=\frac{1}{2^N}$ such that 
\begin{equation}\label{h}h\leq\frac{1}{a{\cal C}}; \end{equation}
it remains to  estimate the complexity of finding an upper bound ${\cal C}$ on $(c_{int}+c_{diff})$ from the input data.

We claim that 
 \begin{equation}\label{c_int} \max\{c_\mathit{int},c_{diff}\}\leq {\cal P}(A,B,C,\varphi),\end{equation}
 where  $${\cal P}(A,B,C,\varphi)=\frac{\lambda_{max}(A)}{\lambda_{min}(A)}\cdot\max\{||\frac{\partial^2
\varphi}{\partial x_i\partial x_j}||_{s}\}\cdot$$
 $$\cdot\max\{||A||_2, ||B||_2,||C||_2,||(A^{-1}B)^2||_2, ||(A^{-1}C)^2||_2,||A^{-1}BA^{-1}C-A^{-1}CA^{-1}B||_2\}.$$

As noted in
\cite{gr} (see chapter 5),
$c_{\operatorname{diff}}=c_1\cdot c_2$ where $c_2=\sqrt{\frac{\lambda_{max}(A)}{\lambda_{min}(A)}}$ comes from the
stability property \cite{go76} and $c_1$ comes from the approximation property
$||L_hu_h-(Lu)|_{G^\tau_k}||_{sL_2}\leq c_1h$. Since our scheme has the
first order of approximation, it follows from the Taylor
decomposition of $L{\bf u}$ that $c_1$ depends only on the norms of $A,B_i$ and
 $$||\frac{\partial {\bf u}}{\partial x_i}||_{sL_2},\; ||\frac{\partial
{\bf u}}{\partial t}||_{sL_2},\; ||\frac{\partial^2 {\bf u}}{\partial
x_i\partial x_j}||_{sL_2},\; ||\frac{\partial^2 {\bf u}}{\partial
x_i\partial t}||_{sL_2}.$$

As it is known, by the proof of the uniqueness theorem for
\eqref{sist_1} \cite{go71}  (p. 155 for the Cauchy problem and p. 194 for the boundary-value problem, respectively), see also \cite{evans,mi}, we have 
 \begin{equation}\label{energy}||{\bf u}||_{A,sL_2}\leq ||\varphi||_{A,L_2}\end{equation} 
   Applying an analogue of \eqref{energy} to the systems for the second derivatives of {\bf u}  and using the equivalence of norms in $\mathbb R^n$ $$\lambda_{min}(v,v)\leq||v||_A=(Av,v)\leq\lambda_{max}(v,v),$$ we obtain the desired estimate.
More precisely, the estimates for the derivatives of the solution can be obtained as follows. 
Considering the Cauchy problem,  due to
the smoothness assumptions, we can construct auxiliary Cauchy
problems for partial derivatives of ${\bf u}$ (we write down a couple of
them, as examples):
 $$
\begin{cases} A({\bf u}_x)_t+B({\bf u}_x)_x+C({\bf u}_x)_y=0,\\
{\bf u}_x|_{t=0}=\varphi_x, \quad
\end{cases}
\begin{cases} A({\bf u}_t)_t+B({\bf u}_t)_x+C({\bf u}_t)_y=0,\\
{\bf u}_t|_{t=0}=-A^{-1}(B\varphi_x+C\varphi_y),
\end{cases}
 $$
 \begin{equation} \label{sist_13}
\begin{cases} A({\bf u}_{tt})_t+B({\bf u}_{tt})_x+C({\bf u}_{tt})_y=0,\\
{\bf u}_{tt}|_{t=0}=-A^{-1}(B({\bf u}_t|_{t=0})_x+C({\bf u}_t|_{t=0})_y)=\\=(A^{-1}B)^2\varphi_{xx}+(A^{-1}C)^2\varphi_{yy}+(A^{-1}BA^{-1}C-A^{-1}CA^{-1}B)\varphi_{xy}.
\end{cases}
\end{equation}
  From \eqref{int-est} $c_{int}$ is easily estimated in a similar way as above (by bounding the derivatives of ${\bf u}$). 
  
  From these considerations, $h\leq\frac{1}{a{\cal C}}\leq\frac{1}{2a{\cal P}(A,B,C,\varphi)}$ can be computed in PTIME, since all the expressions in ${\cal P}(A,B,C,\varphi)$ (eigenvalues, matrix multiplication, taking an inverse matrix, calculating upper bounds on the norms, differentiating rational polynomials) are PTIME computable.
   $\qed$

Note that for not fixed  $n$, the PTIME bounds in Proposition \ref{steps_prop} fail to hold, since there are no PTIME bounds on finding eigenvalues, see Remark \ref{rem_spec}.

\subsection{Description of the algorithm}

The difference scheme for the boundary-value problem
\eqref{sist_2} and the Cauchy problem  \eqref{sist_1} may be
chosen in various ways. We use the  Godunov scheme \cite{go76} (see also e.g. \cite{kps}), which 
can also be applied to a broader class of systems, including  some systems of
nonlinear equations. We describe it in few stages, letting for
simplicity the righthand part of \eqref{sist_1} to be zero: $f=0$. 
 The scheme approximates the system \eqref{sist_1} or
\eqref{sist_2} with the first order of accuracy  (the
proof of the approximation property is done by means of the Taylor
decomposition). To simplify the analysis of complexity, we write the scheme in an algorithmic form. 

The algorithm describes, in several steps, computation of the values ${\bf
u}^{i-\frac{1}{2},j-\frac{1}{2}}$ (approximating the solution of \eqref{sist_1} or \eqref{sist_2} at the point $(t,\frac{i-\frac{1}{2}}{2^N},\frac{j-\frac{1}{2}}{2^N})$) on the time level $t=(l+1)\tau$ from the values ${\bf
u}_{i-\frac{1}{2},j-\frac{1}{2}}$ on the time level $t=l\tau$. We distinguish  the values by upper and lower indices, in order to avoid using of a third index $l$. 

0. Find the matrices $T_x$, $T_y$, $T_x^{-1}$, $T_y^{-1}$, $A^{-1}$, the eigenvalues $\mu_k(B)$, $\mu_k(C)$ of the matrix pencils $\mu A-B$, $\mu A-C$ ($k=1,2,\ldots,N$) and the steps $h,\tau$ from  given $m,n,a$, algebraic  (resp. rational) matrices $A, B, C$ and the initial function $\varphi$. 

1. Let (for $i$ from $1$ to  $2^N$, $j$ from $1$ to $2^N$) ${\bf
u}_{i-\frac{1}{2},j-\frac{1}{2}}=\varphi\mid_{G_N}$, i.e. on the level t=0 we just take the values of initial conditions $\varphi(\frac{i-\frac{1}{2}}{2^N},\frac{j-\frac{1}{2}}{2^N})$.

The next steps 2-5 are carried out for all $l$ from $1$ to $L=\frac{T}{\tau}$.

2.  For $i$ from $1$ to  $2^N$, $j$ from $1$ to $2^N$, find the auxiliary vectors  $${\bf v}^{(x)}_{i-\frac{1}{2}, j-\frac{1}{2}}=T_x^{-1}{\bf u}_{i-\frac{1}{2}, j-\frac{1}{2}},\quad {\bf v}^{(x)}_{i+\frac{1}{2}, j-\frac{1}{2}}=T_x^{-1}{\bf u}_{i+\frac{1}{2}, j-\frac{1}{2}}$$
$${\bf v}^{(y)}_{i-\frac{1}{2},j-\frac{1}{2}}=T_y^{-1}{\bf u}_{i-\frac{1}{2}, j-\frac{1}{2}},\quad {\bf v}^{(y)}_{i-\frac{1}{2}, j+\frac{1}{2}}=T_y^{-1}{\bf u}_{i-\frac{1}{2}, j+\frac{1}{2}}.$$

3. Find the auxiliary ``large values'' $${\mathcal V}^{(x)}_{i,j-\frac{1}{2}}=\left(\begin{array}{c}{\mathcal W}_{i,j-\frac{1}{2}}^1\\ {\mathcal W}_{i,j-\frac{1}{2}}^2\\ \ldots \\ {\mathcal W}_{i,j-\frac{1}{2}}^n\end{array}\right);\quad
{\mathcal V}^{(y)}_{i-\frac{1}{2},j}=\left( \begin{array}{c}{\mathcal W}_{i-\frac{1}{2},j}^1\\ {\mathcal W}_{i-\frac{1}{2},j}^2\\ \ldots \\ {\mathcal W}_{i-\frac{1}{2},j}^n\end{array}\right)$$
for $i=0,1,2,\ldots,2^N$, $j=1,2,\ldots,2^N$; $i=1,2,\ldots,2^N$, $j=0,1,2,\ldots,2^N$.

For  auxiliary ``interior'' points $i$ from $1$ to $2^N-1$, j from $1$ to $2^N-1$ let
\begin{equation}\label{flux}{\cal W}_{i,j-\frac{1}{2}}^k=
 \begin{cases}({\bf v}^{(x)}_{i-\frac{1}{2},j-\frac{1}{2}})_k,\text{
if }\mu_k(B)\geq0,\\ ({\bf v}^{(x)}_{i+\frac{1}{2},j-\frac{1}{2}})_k,\text{ if }\mu_k(B)<0;\end{cases} {\cal W}_{i-\frac{1}{2},j}^k=
 \begin{cases}({\bf v}^{(y)}_{i-\frac{1}{2},j-\frac{1}{2}})_k,\text{
if }\mu_k(C)\geq0,\\ ({\bf v}^{(y)}_{i-\frac{1}{2},j+\frac{1}{2}})_k,\text{ if }\mu_k(C)<0,\end{cases}
\end{equation}
where $\mu_k(B), \mu_k(C)$ are respectively the $k$-th eigenvalues of the matrix pencils $\mu A-B$, $\mu A-C$.

3a). In the case of the Cauchy problem \eqref{sist_1}, for the auxiliary ``boundary'' values ${\mathcal W}_{0,j-\frac{1}{2}},$ ${\mathcal W}_{2^N,j-\frac{1}{2}}$ (for $j$ from $1$ to $2^N$) and ${\mathcal W}_{i-\frac{1}{2},0},$ ${\mathcal W}_{i-\frac{1}{2},2^N}$ (for $i$ from $1$ to $2^N$)   we use the same formula \eqref{flux} taking 
$${\bf v}^{(x)}_{-\frac{1}{2},j-\frac{1}{2}}={\bf v}^{(x)}_{2^N-\frac{1}{2},j-\frac{1}{2}},\quad  {\bf v}^{(x)}_{2^N+\frac{1}{2},j-\frac{1}{2}}={\bf v}^{(x)}_{\frac{1}{2},j-\frac{1}{2}}$$ and $${\bf v}^{(y)}_{i-\frac{1}{2},-\frac{1}{2}}={\bf v}^{(y)}_{i-\frac{1}{2},2^N-\frac{1}{2}}, \quad {\bf v}^{(y)}_{i-\frac{1}{2},2^N+\frac{1}{2}}={\bf v}^{(y)}_{i-\frac{1}{2},\frac{1}{2}}.$$   

3b). In the case of the boundary-value problem \eqref{sist_2}, we
compute the boundary values $\mathcal{V}_0,\mathcal{V}_{2^N}$
with the help of the boundary conditions. On the left boundary
$x=0$  we calculate $m_{+}$ components of $\mathcal{V}_0$,
corresponding to the positive eigenvalues of the matrix $A^{-1}B$,
from the system of linear equations $\Phi_1^{(1)}\mathcal{V}_0=0$;
for $m_{-}$ components of $\mathcal{V}_0$, corresponding to the
negative eigenvalues, we let $\mathcal{V}_0:=v_{\frac{1}{2}}$. The
components corresponding to the zero eigenvalues of $A^{-1}B$ can
be chosen arbitrarily since they are multiplied by zero in the
scheme. The values on the right boundary and on both boundaries by
the $y$-coordinate are calculated in a similar way.

4. Calculate ${\cal U}_{i,j-\frac{1}{2}}=T_x{\cal V}^{(x)}_{i,j-\frac{1}{2}}$ (for $i$ from $0$ to $2^N$, for $j$ from $1$ to $2^N$) and ${\cal U}_{i-\frac{1}{2},j}=T_y{\mathcal V}^{(y)}_{i-\frac{1}{2},j}$ (for $i$ from $1$ to $2^N$, for $j$ from $0$ to $2^N$).

5. Find values on the next grid step: for $i$ from $1$ to $2^N$, $j$ from $1$ to $2^N$  let 

\begin{equation}\label{scheme_2_sist1} {\bf
u}^{i-\frac{1}{2},j-\frac{1}{2}}={\bf
u}_{i-\frac{1}{2},j-\frac{1}{2}}-\frac{\tau}{h}A^{-1}\left(B({\cal
U}_{i,j-\frac{1}{2}}-{\cal
U}_{i-1,j-\frac{1}{2}})+C({\cal
U}_{i-\frac{1}{2},j}-{\cal U}_{i-\frac{1}{2},j-1})\right).
\end{equation}

Remember the calculated values as ${\bf
u}_{i-\frac{1}{2},j-\frac{1}{2}}^{(l)}={\bf
u}^{i-\frac{1}{2},j-\frac{1}{2}}$, then let ${\bf
u}_{i-\frac{1}{2},j-\frac{1}{2}}={\bf
u}^{i-\frac{1}{2},j-\frac{1}{2}}$.

6. Finally, $\upsilon=\left\{{\bf
u}_{i-\frac{1}{2},j-\frac{1}{2}}^{(l)}\right\}_{l=1}^L\mid_{H}$ is the approximation of the solution ${\bf u}$ of the system \eqref{sist_1}, and $\upsilon=\left\{{\bf
u}_{i-\frac{1}{2},j-\frac{1}{2}}^{(l)}\right\}_{l=1}^L\mid_{[0,1]^m\times[0,T]}$ is the approximation of the solution ${\bf u}$ of the system \eqref{sist_2}.

\begin{Remark} Stage 3 of the algorithm can be carried out without using the branching operator, by letting 
$${\cal W}_{i,j-\frac{1}{2}}=
 S^B_{-}{\bf v}^{(x)}_{i+\frac{1}{2},j-\frac{1}{2}}+S^B_{+}{\bf v}^{(x)}_{i-\frac{1}{2},j-\frac{1}{2}};\ {\cal W}_{i-\frac{1}{2},j}=
 S^C_{-}{\bf v}^{(y)}_{i-\frac{1}{2},j+\frac{1}{2}}+S^C_{+}{\bf v}^{(y)}_{i-\frac{1}{2},j-\frac{1}{2}}.$$ 
Here $S^B_{-}$ (resp. $S^C_{-}$) is the matrix $\text{diag}\{1,1,\ldots,1,0,0,\ldots,0\}$,  with the number of $1$s  equal to the number of negative eigenvalues of $A^{-1}B$ (resp. $A^{-1}C$); similarly, $S^B_{+}$ (resp. $S^C_{+}$) is the matrix $\text{diag}\{0,0,\ldots,0,1,1,\ldots,1\}$ with the number of $1$s  equal to the number of nonnegative eigenvalues of $A^{-1}B$ (resp. $A^{-1}C$). Note that the matrices $S^B_{-}$, $S^C_{-}$, $S^B_{+}$, $S^C_{+}$ can be computed on Stage 0 before the cycles, and that these matrices depend only on the signs of the eigenvalues $\mu_k(B)$, $\mu_k(C)$.
 \end{Remark}

\subsection{Counting steps}\label{csteps}

With all this at hand, it is not hard to count the computation steps in the Godunov scheme (all computations are w.r.t. the $p$-presentation of $\mathbb{A}$ in Section \ref{algprel}). By Proposition \ref{steps_prop}, the number of grid points in the scheme (see Section \ref{discr}) is bounded by a polynomial. The computations in Godunov's scheme proceed bottom-up by layers, along the time axis. At the bottom level, we just evaluate the initial functions in the grid points which requires polynomial time according to remarks in Section \ref{algprel}. To go one level up requires, for each grid point on the next level, the values at the previous levels and a fixed number of matrix multiplications by matrices, computed in advance using Propositions \ref{specmat} and \ref{specpen}. Therefore, climbing one level up also requires polynomial time. Let $p_i$, $i=1,\ldots,L$ (where $L=\frac{T}{\tau}$, computed in polynomial time), be a polynomial bounding the computation time for level $i$. Since the computation at level $i$ uses only the values of $\upsilon$ at grid points of level $i-1$ (note that computation of the value at any point from the $i$-th level requires only finite number of points at the ($i-1$)th level, in our case five $2m+1$ points) and some matrices computed in advance, the whole computation time is (essentially) bounded by the composition $p_L\circ\cdots\circ p_1$ of polynomials which lays down to EXPTIME and yields the complexity bound given in item 1 of Theorem \ref {cpq}.

For item 2 the argument is the same, except that one has to take into account  that for a fixed $n$ the algorithm of spectral decomposition works in polynomial time and that the number of time steps $L$ is just constant (instead of PTIME computable)  in this case.  
It follows from the estimates \eqref{tau}, \eqref{h} and \eqref{c_int} of Proposition \ref{steps_prop}, in particular the expression ${\cal P}(A,B,C,\varphi)$ used to calculate $h$ can be taken just as $M^3$.
$\qed$
%For item 3  we have additionally to use the corresponding items of  Proposition \ref{specmat}. This completes the proof of Theorem \ref {cpq}. 

\begin{Remark} 1) Note that the estimate in Theorem \ref {cpq} (1) is exponential even for  fixed $n$  because we take arbitrary rational polynomials as initial functions. Taking reasonably restricted classes of initial functions and matrices (with restrictions like those in \cite{ss09,ss13}) yields a polynomial estimate in item (2), though the degree of polynomial is high and one needs to take a fine grid with small steps $h$ and $\tau$. In this way, the exponential algorithm of item (1) might work out better than the polynomial one of item (2) for concrete problems. Also note that in item (2) the precision is fixed.

2) For the case of simultaneously diagonalizable (or, equivalently, mutually commuting) matrices $A^{-1}B_j$, the solution of the IVP \eqref{sist_1} even in the sense of Task 1 is in PTIME (without using the difference scheme method). Indeed, for $n=1$, the scalar transport equation 
$u_t=\sum\limits_{j=1}^m b_ju_{x_j}$ can be easily solved in PTIME: $u(t,x_1,\ldots,x_m)=\varphi(x_1-b_1t,
\ldots,x_m-b_mt)$.  For a fixed $n>1$, the system \eqref{sist_1} can be linearly transformed to $n$ independent transport equations via spectral decomposition of the matrix pencils $(A,B_j)$, which is computable in PTIME according to Proposition \ref{specpen}.
   \end{Remark}

\section{Conclusion}\label{con}

In this paper we obtained apparently first  bit complexity upper bounds for  computing solutions of the Cauchy and dissipative boundary-value problems for symmetric hyperbolic systems of PDEs (to which also many higher-order hyperbolic PDEs can be reduced) with guaranteed precision. 

Although our methods do not always yield (for instance, for large $n$) practically feasible algorithms for guaranteed-precision problems for PDEs, we  hope  that investigations in this direction are fruitful for both theoretical research and applications. In particular, on the implementation level it seems useful and rewarding to enhance the existing systems of ``exact real computations'' by packages based on highly developed algorithms of computer algebra. We are not aware of the existence of such ``hybrid'' systems built under the slogan of ``guaranteed-precision numerical computations''.

{\bf Acknowledgement.} We are grateful to Sergey  Goncharov for stimulating discussions, and to Pavel Alaev, Gennadiy Demidenko and Martin Ziegler for valuable comments.


\begin{thebibliography}{aasaaa67}

\bibitem{ak89} A.G.Akritas,  \textit{Elements of Computer Algebra with Applications}.  Wiley Interscience, New York, 1989.

\bibitem{al16} P.E. Alaev, Existence and uniqueness of structures computable in polynomial time, \textit{Algebra and Logic}, \textbf{55}, No 1 (2016),  106--112. 
% \bid{doi={10.1007/s10469-016-9377-6}}

\bibitem{as17} P.E. Alaev, V.L. Selivanov, Fields of algebraic numbers computable in polynomial time. I, Algebra and Logic, 58:6 (2019), 673--705.

%P.E. Alaev  and V.L. Selivanov,  \textit{Polynomial-time presentations of algebraic number fields (extended abstract)}, in: LNCS volume  10936 of Proceedings of the conference Computability in Europe (Ed. F. Manea, R. Miller and . Novotka), 2018, Berlin, Springer, pp. 20--29.  \bid{doi={10.1007/978-3-319-94418-02}}


\bibitem{ba86} K.I. Babenko,    \textit{Foundations of Numerical Analysis}.
Moscow, Nauka, 1986. (in Russian).

\bibitem{ber} P. B\"urgisser, M. Clausen, A. Shokrollahi, \textit{Algebraic Complexity Theory}. Berlin-Heidelberg, Springer, 1997.

\bibitem{bdg} J.L. Balc\'{a}zar, J. D\'{i}az, and J. Gabarr\'{o}.
\textit{ Structural Complexity I}, in: volume 11 of EATCS Monographs on
Theoretical Computer Science, Springer-Verlag, 1988.

\bibitem{bhw}  V. Brattka,  P. Hertling and K. Weihrauch,
\textit{A tutorial on computable analysis}, in: New Computational Paradigms
(edited by S. Barry Cooper, Benedikt L\"owe, Andrea Sorbi), 2008,
pp. 425--491.

\bibitem{bp06}  S. Basu, R. Pollack and M. Roy, 
\textit{Algorithms in Real Algebraic Geometry}, Springer, Heidelberg, 2006.

\bibitem{cl82}  G.E. Collins  and R. Loos, \textit{ Real zeros of polynomials}, in: ``Computer Algebra: Symbolic and Algebraic Computations'', Springer-Verlag, 1982,  pp. 83--94.


\bibitem{cr91}D.  Cenzer and J. Remmel, Polynomial time versus recursive models, \textit{ Annals of Pure and Applied Logic},\textbf{  54} (1991), 17--58.
%\bid{doi={10.1016/0168-0072(91)90008-a}}

\bibitem{eg99}  Yu.L. Ershov  and  S.S. Goncharov, \textit{ Constructive
Models}, Novosibirsk, Scientific Book, 1999 (in Russian, there is an
English Translation).

\bibitem{evans} L.C. Evans, \textit{ Partial Differential
Equations}, in: Graduate Studies in Mathematics, v. 19,  American
Mathematical Society, 1998.

\bibitem{fri} K.O. Friedrichs,  Symmetric hyperbolic linear differential
equations, \textit{ Communication on Pure and Applied Mathematics}, \textbf{ 7}
(1954), 345--392.

\bibitem{gant} F.R. Gantmacher, \textit{  Matrix Theory},
Nauka, Moscow, 1967 (in Russian).

\bibitem{gr} S.K. Godunov  and V.S. Ryaben'kii, \textit{ ntroduction
to the
Theory of Difference Schemes}. Fizmatgiz, Moscow, 1962 (in Russian). English translation:
\textit{ Difference Schemes: An Introduction to the  Underlying Theory
(Studies in Mathematics and Its Applications)}, Elsevier Science
Ltd, 1987.


\bibitem{God} S.K. Godunov  et al., \textit{ Guaranteed Precision of Solving Systems of Linear Equations in Euclidean Spaces}, Novosibirsk, Nauka, 1988 (in Russian).

\bibitem{go71} S.K. Godunov, \textit{ Equations of Mathematical Physics}, Nauka, Moscow, 1971 (in Russian).

\bibitem{go76} S.K. Godunov, ed. \textit{ Numerical Solution of
Higher-dimensional Problems of Gas Dynamics}  Nauka,
Moscow, 1976 (in Russian).

\bibitem{hj} R.A. Horn  and Ch. R. Johnson, \textit{ Matrix analysis},
Cambridge University Press, 1983.

\bibitem{jo66} F. John, \textit{  Lectures on Advanced Numerical Analysis},
Gordon and Breach, Science Publishers, Inc., 1966.

\bibitem{ko} Ko Ker-I, \textit{  Complexity Theory of Real Functions},
Birkh\"auser, Boston, 1991.

\bibitem{kz} A. Kawamura, M. Ziegler, \textit{ Invitation to Real Complexity Theory: Algorithmic Foundations to Reliable Numerics with Bit-Costs}, arXiv:1801.07108.

\bibitem{ksz}  I. Koswara,  S. Selivanova and M. Ziegler, \textit{  Computational complexity of real powering and improved solving linear differential equations}, in: Proc. 14th InternationalComputer Science Symposium in Russia, vol. 11532 of LNCS, 2019.

\bibitem{kps} A.G.  Kulikovskii, N.V. Pogorelov and A.Yu. Sem\"enov,
 \textit{ Mathematical Aspects of Numerical Solution of Hyperbolic
Systems}, Chapman \& Hall/CRC Press, Boca Raton, 2001.


\bibitem{lll82} A.K. Lenstra, H.W. Lenstra and L. Lovasz, Factoring polynomials with rational coefficients,  \textit{ Math. Ann.},  \textbf{ 261} (1982),  515--534. %\bid{doi={10.1007/BF01457454}}



\bibitem{lo82} R. Loos,  \textit{  Computing in algebraic extensions}, in: ``Computer Algebra: Symbolic and Algebraic Computations'', Springer-Verlag, 1982, pp. 115--138.


\bibitem{mi} S. Mizohata,  \textit{  The Theory of Partial Differential
Equations}, Cambridge Univ. Press, Cambridge, 1973.

\bibitem{panchen} V.Y. Pan and Z.Q. Chen,  \textit{ The complexity of the matrix eigenproblem}, in: Proceedings of the thirty-first annual ACM symposium on Theory of computing, ACM, 1999, pp. 507--516.

\bibitem{pan} V. Pan  and J. Reif,  The bit complexity of discrete solutions of partial differential equations: compact multigrid,  \textit{ Computers Math. Applic.},  \textbf{ 20}, No 2 (1990), 9--16. %\bid{doi={10.1016/0898-1221(90)90235-C}}


\bibitem{rel}  F. Rellich,  St\"orungstheorie der Spektralzerlegung I., Analytische Sto\"rung der
isolierten Punkteigenwerte eines beschr\"ankten Operators,  \textit{  Math. Ann.},  \textbf{ 113}, 1937, 600--619 (in German).

\bibitem{sc86} A. Schrijver,  \textit{Theory of Linear and Integer Programming},  Wiley Interscience, New York, 1986.

\bibitem{ss09}  S. V. Selivanova and V. L. Selivanov,
Computing solution operators of symmetric hyperbolic systems of
PDEs, \textit{Journal of Universal Computer Science}, \textbf{15:6} (2009),
 1337--1364.

\bibitem{ss13} S.  Selivanova  and V. Selivanov, Computing solution operators of boundary-value
problems for some linear hyperbolic systems of PDEs, \textit{Logical Methods in Computer Science}, \textbf{4:13} (2017), 1--31. Earlier version on  arXiv:1305.2494 (2013).

\bibitem{ss17} S.V. Selivanova  and V.L. Selivanov, On constructive number fields and
computability of solutions of PDEs, \textit{Doklady Mathematics},  \textbf{477:3} (2017), 282--285. 

\bibitem{ss18} S.V. Selivanova and  V.L. Selivanov, \textit{Bit complexity  of computing  solutions  for symmetric hyperbolic systems of PDEs  
(extended abstract)}, in: LNCS volume  10936 of Proceedings of the conference Computability in Europe (Ed. F. Manea, R. Miller and D. Novotka),  Berlin, Springer, 2018.

\bibitem{sz}
S. Selivanova  and M. Ziegler, \textit{Turnkey Solutions to PDEs in Exact Real Computation} in: Book of Abstracts - 12th Summer Workshop on Interval methods, Palaiseau, France, July 23-26, 2019, pp. 21--22. https://swim2019.ensta-paris.fr/SWIM-19$\_$Book$\_$of$\_$Abstracts.pdf

\bibitem{sma} S. Smale, On the topology of algorithms, \textit{ J. of Complexity},  \textbf{3} (1987), 81--89.

\bibitem{strik} J.C. Strikverda, \textit{Finite Difference Schemes and
Partial Differential Equations}, SIAM, 2004.


\bibitem{sz59} G. Szeg\"o, \textit{Orthogonal Polynomials},
AMS 531, West 116th Street, New York, 1959.


\bibitem{tre96} L.N.  Trefethen, \textit{Finite Difference and
Spectral Methods for Ordinary and Partial Differential Equations},
Cornell University, Department of Computer Science and Center for
Applied Mathematics, 1996.

\bibitem{vas} V. Vassiliev, Cohomology of braid groups and the complexity of algorithms, \textit{ Functional Anal. Appl}, \textbf{22} (1989), 182--190.

\bibitem{wa67}  B.L. van der Waerden, \textit{ Algebra},
Springer, Berlin, 1967.

\bibitem{wei} K. Weihrauch, \textit{ Computable Analysis}, Berlin,
Springer, 2000.

\bibitem{wei03} K. Weihrauch, Computational Complexity on Computable Metric Spaces, \textit{ Mathematical Logic Quarterly}, \textbf{49:1} (2003), 3--21. 


\bibitem{yap} C. Yap, M. Sagraloff, V. Sharma, \textit{Analytic root clustering: A complete algorithm using soft zero tests}, in: Conference on Computability in
Europe, Springer, 2013, pp. 434--444.

\bibitem{zhou} J.-P. Zhou, On the degree of extensions generated by finitely many algebraic numbers \textit{ Journal of Number Theory}, \textbf{34} (1990), 133--141.

\bibitem{za} Yu.S. Zavyalov, B.I. Kvasov and V.L. Miroshnichenko,
\textit{ Methods of the Spline Functions}, Fizmatgiz, Moscow,
1980 (in Russian).

\bibitem{zie}  M. Ziegler,  Real computation with least discrete advice: a complexity theory of nonuniform
computability, \textit{ Annals of Pure and Applied Logic}, \textbf{163:8} (2012), 1108--1139.

\bibitem{zb04} M. Ziegler and V. Brattka, \textit{ A computable spectral
theorem}, in:
Proc. CCA-2001, Lecture Notes in Computer Science, v. 2064
(2001), 378--388.


\end{thebibliography}
\end{document}